\documentclass[10pt, conference]{ieeeconf}
\IEEEoverridecommandlockouts                              % This command is only needed if 
\usepackage{verbatim}
\usepackage{float}
\usepackage{mathrsfs}
\usepackage{amsmath}
\usepackage{amssymb}

\makeatletter

\floatstyle{ruled}
\newfloat{algorithm}{tbp}{loa}
\providecommand{\algorithmname}{Algorithm}
\floatname{algorithm}{\protect\algorithmname}

\newtheorem{thm}{Theorem}
\newtheorem{problem}{Problem}
\newtheorem{rem}{Remark}
\newtheorem{defn}{Definition}
\newtheorem{lem}{Lemma}
\newtheorem{cor}{Corollary}

\makeatother

\usepackage{babel}

\title{Nonlinear Dynamical Unbalanced Optimal Transport: Relaxation and Duality}
\author{Dongjun Wu, Anders Rantzer
	\thanks{*This project has received funding from the European Research Council (ERC) under the European Union's Horizon 2020 research and innovation programme under grant agreement No 834142 (ScalableControl).}
	\thanks{The authors are with the Department of Automatic Control, Lund
	University, Box 118, SE-221 00 Lund, Sweden {\tt\small
	(dongjun.wu, anders.rantzer)@control.lth.se}.}
}

\begin{document}
\maketitle
\begin{abstract}
In this paper, we introduce a generalized dynamical unbalanced optimal
transport framework by incorporating limited control input and mass
dissipation, addressing limitations in conventional optimal transport
for control applications. We derive a convex dual of the problem using
dual optimal control techniques developed before and during the 1990s,
transforming the non-convex optimization into a more tractable form.
At the core of this formulation is the smooth sub-solutions to an HJB
equation. A first-order algorithm based on the dual formulation is
proposed to solve the problem numerically. 
\end{abstract}
{\keywords Nonlinear systems, Optimal transport, dual optimal control}

\section{Introduction} \label{sec:intro}

Optimal transport (OT) has emerged as a powerful mathematical framework
with growing interest in the control community \cite{chenOptimalTransportSystems2021,bunne2023optimal},
particularly for its ability to model and optimize the transportation
of resources, distributions, or probabilities in a variety of applications
such as robotic swarms \cite{elamvazhuthiMeanfieldModelsSwarm2019},
motion planning \cite{le2023accelerating} and single-cell analysis
\cite{bunne2023learning}. 

Despite its successes, conventional OT, rooted in Monge\textquoteright s
formulation and Kantorovich\textquoteright s relaxation, assumes mass
conservation, meaning the total mass of the initial and target measures
must be equal. This assumption limits its applicability in scenarios
where mass changes occur during the transport process. For example,
in the transmission of electricity through a grid, energy losses due
to resistance lead to mass dissipation; in irrigation systems, evaporation
or leakage causes mass reduction; and in biological systems, cell
growth or decay introduces mass creation or destruction. These cases
highlight a critical limitation of classical OT: its inability to
handle mass variations, which are prevalent in control applications.
To overcome this limitation, extensions to incorporate mass creation
and dissipation have been proposed in the literature \cite{chizatInterpolatingDistanceOptimal2018,figalliOptimalPartialTransport2010,Kondratyev2016,lieroOptimalTransportCompetition2016,liero2018optimal}.
We adopt the term ``unbalanced optimal transport'' (UOT) for such
extended OT problems. UOT often uses partial differential equations
(PDEs) to govern the evolution of measures with varying total mass,
as seen in the dynamical formulation of Hellinger-Kantorovich distance
(HK) \cite{liero2018optimal} or Wasserstein-Fisher-Rao (WF) metric
\cite{chizatInterpolatingDistanceOptimal2018}. While UOT has expanded
the scope of OT to dynamical systems with mass variation, existing
formulations typically assume unconstrained control inputs and mass
dissipation, which may not reflect practical constraints (both dynamical
and static) in control problems. For instance, in robotic swarms,
nonlinear dynamics are present and control inputs (e.g., torques and
velocities) are often bounded by physical limits, and in biological
systems, mass changes (e.g., growth rates) are constrained by environmental
or biochemical factors.

In this paper, we propose a generalized formulation of dynamical UOT
that explicitly accounts for systems dynamics, constrained control
inputs and constrained mass changes. The primary goal of this work
is to derive a convex duality for this generalized UOT problem, transforming
the nonlinear optimization into a dual form that is more tractable
for analysis and computation. To achieve this, we leverage techniques
from classical dual optimal control theory, that have been developed
before and during the 1990s \cite{young1969lectures,flemingConvexDualityApproach1989,vinterConvexDualityNonlinear1993}.
Specifically, we employ the concept of relaxed control to lift the
problem into a convex setting. This approach not only establishes
a duality result but also enables the development of a first-order
numerical algorithm to solve the problem. 

Our contributions provide a flexible framework for control problems
involving mass variation, with potential applications in energy systems,
fluid dynamics, and biological modeling.

{\it Notations:} We collect some notations that will be used in the paper. Throughout the paper, $X$ is a compact set in $\mathbb{R}^n$. $\partial X$ denotes the boundary of $X$. The set of non-negative Borel measures on $X$ is denoted $\mathcal{M}_+(X)$; the set of probability measures on $X$ is denoted $\mathscr{P}(X)$.
Given a function $f:[0,1]\times X \to \mathbb{R}$, $\nabla f$ typically means the gradient with respect to the second variable.
Let $X_1$ and $X_2$ be two measurable spaces, $\mu$ a measure on $X_1$, and $T$ a measurable map from $X_1$ to $X_2$. Then $T_{\#} \mu$ is the push-forward measure on $X_2$. Given a set $A \subseteq \mathbb{R}^n$, the distance from a point to the set $A$ is defined as $d_A(x)  = \inf_{y\in A} |x-y|$. Given a normed linear space $\mathcal{X}$, the set of bounded linear functionals on $\mathcal{X}$ is denoted $\mathcal{X}^*$. Given a set $S \subseteq \mathcal{X}$, ${\rm co} S$ is the smallest convex set in $\mathcal{X}$ containing $S$, and $\overline{\rm co}(S)$ is the closure of ${\rm co}(X)$.

\section{Problem setting}

Let $X$ be a compact set in $\mathbb{R}^{n}$. Consider the following
system described by a partial differential equation (PDE) in variable
$\mu(t,x)$:
\begin{equation}
\partial_{t}\mu+\nabla\cdot(\mu f(t,x,u)=g(t,x,u)\mu\label{sys:fg}
\end{equation}
with boundary condition
\begin{equation}
\mu(0,x)=\mu_{0}(x),\quad\mu(1,x)=\mu_{1}(x)\label{eq:bound}
\end{equation}
The unknown variable $\mu(t,\cdot)$ is a one-parameter non-negative
Borel measure, i.e., $\mu(t,\cdot)\in\mathcal{M}_{+}(X)$ and is compactly
supported in $X$ for every $t\in[0,1]$. The vector field $f(t,x,\cdot)$
and function $g(t,x,\cdot)$ are given and $u(t,x)$ is to be designed
to achieve certain optimal criteria. We call $u$ the control input.
For convenience, denote $D=[0,1]\times X$.

The equation (\ref{sys:fg}) should be understood in distributional
sense. Namely, the following identity holds for all $\varphi\in C^{1}([0,1]\times X)${\small
\[
\int_{0}^{1}\int_{X}(\partial_{t}\varphi+\nabla\varphi\cdot f+g\varphi){\rm d}\mu(t,x){\rm d}t=\int_{X}\varphi(1,\cdot){\rm d}\mu_{1}-\varphi(0,\cdot){\rm d}\mu_{0}.
\]
We shall call (\ref{sys:fg}) the }{\small\emph{system}}{\small , which
}describes how an initial measure $\mu_{0}$ is transported to $\mu_{1}$
under the flow of the system
\begin{equation}
\dot{x}=f(t,x,u(t,x))\label{sys:ode}
\end{equation}
with mass dissipation determined by $g(\cdot,\cdot,\cdot)$. Indeed,
in case $\mu(t,x)$ is a compactly supported density on $X$, then
the mass evolution is governed by
\begin{align*}
\frac{d}{dt}\int_{X}\mu(t,x){\rm d}x & =\int_{X}-\nabla\cdot(\mu f)+g\mu{\rm d}x=\int_{X}g\mu{\rm d}x.
\end{align*}
If $g\equiv0$, then the total mass $\int\mu(t,\cdot)$ is preserved,
and $\mu(t,x)$ is the pushforward measure of $\mu_{0}$ by the flow
of the system (\ref{sys:ode}). If $g\ge0$ (resp. $\le0$), then
the total mass is monotone increasing (resp. decreasing). 

With (\ref{sys:fg}) being our model and (\ref{eq:bound}) the boundary
condition, we study an optimization problem of the following type:
\begin{equation}
\min_{u,\mu}\int_{0}^{1}\int_{X}L(t,x,u){\rm d}\mu(t,x){\rm d}t\label{eq:cost}
\end{equation}
where $L(\cdot,\cdot,\cdot)$ is a non-negative function called the
Lagrangian. 

\begin{comment}
with
\[
v(t,x)\in\mathcal{F}(t,x)\subseteq\mathcal{F},\quad w(t,x)\in\mathcal{G}(t,x)\subseteq\mathcal{G}
\]
for some compact convex sets $\mathcal{F}(t,x),\mathcal{G}(t,x)$
and $\mathcal{F},\mathcal{G}$. 
\end{comment}
\begin{comment}
We describe a typical example of the preceding sets. 
\end{comment}
\begin{comment}
The most important case: $f(t,x,u)=f_{0}(t,x)+f_{1}(t,x)u$ and $g(t,x,u)=g_{0}(t,x)+g_{1}(t,x)u$
where $f_{1}$ (resp. $g_{1}$) is a matrix (resp. row vector), and
that $u\in U\subseteq\mathbb{R}^{m}$ for some compact set $U$. Assume
that $f_{i},g_{i}$ are all locally Lipschitz in $x$ and piece-wise
continuous in $t$. In this paper, we consider measurable input $u$
such that the following regularity assumption is satisfied.
\end{comment}

Under certain regularity assumption, the system (\ref{sys:fg}) is
known to have a unique weak solution in the distributional sense given
initial datum $\mu_{0}\in\mathscr{P}(X)$. And the solution is fully
characterized by the flow of the ordinary differential equation (\ref{sys:ode})
and the mass dissipation term $g$: for given $u$ with required regularity,
let $X_{u}(t;x)$ be the system flow of (\ref{sys:ode}) starting
from $x$ at $t=0$, then {\small
\begin{equation}
\mu(t,\cdot)=X_{u}(t,\cdot)_{\#}\left(\mu_{0}\exp\left(\int_{0}^{t}g(s,X_{u}(s,\cdot),u(t,X_{u}(s,\cdot)))\right)\right)\label{eq:sol mu(t,)}
\end{equation}
}for $t\in[0,1]$, see for example \cite{manigliaProbabilisticRepresentationUniqueness2007}.
In particular, if $\mu_{0}=\delta_{x_{0}}$, i.e., a Dirac mass at
$x_{0}$. Then the solution to (\ref{sys:fg}) is
\[
\mu(t,x)=e^{\int_{0}^{t}g(s,x(s))ds}\delta_{x(t)}(x)
\]
where $x(t)$ is the solution to (\ref{sys:ode}) with initial condition
$x(0)=x_{0}$.

We ignore the technical regularity issues concerning the well-posedness
of the system solution, and instead assume that the solution to the
system always exists and is unique once the input and initial data
$\mu_{0}$ have fixed to certain set (to be defined in the sequel).

To facilitate the study of the problem (\ref{eq:cost}) under the
dynamic constraints (\ref{sys:fg}), we reformulate it into the following
form.
\begin{problem}[Primal problem]
\label{pb:P}Consider the system
\begin{equation}
\begin{aligned}\partial_{t}\mu(t,x)+\nabla\cdot[\mu(t,x)v(t,x)] & =w(t,x)\mu(t,x)\\
\mu(0,x) & =\mu_{0}(x)\\
\mu(1,x) & =\mu_{1}(x)
\end{aligned}
\label{sys:mu_v_w}
\end{equation}
Study the Lagrangian type problem
\begin{equation}
\min_{v,w,\mu}\int_{0}^{1}\int_{X}L(t,x,v(t,x),w(t,x)){\rm d}\mu(t,x){\rm d}t\label{eq:cost-Lag}
\end{equation}
subject to constraints
\begin{equation}
\begin{bmatrix}v(t,x)\\
w(t,x)
\end{bmatrix}\in F(t,x)\label{eq:sys-cnst}
\end{equation}
where the function $L$ and the set $F$ satisfy Assumption H1 below.
\end{problem}
%
\begin{comment}
For example, if $f(t,x,u)=f(t,x)+h(t,x)u$, $g(t,x,u)=a(t,x)+b(t,x)u$,
and that $u(t,x)\in U$, then 
\end{comment}

\textbf{Assumption H1}: 
\begin{itemize}
\item For every $(t,x)\in D$, $F(t,x)$ is a convex set, and $(v,w)\mapsto L(t,x,v,w)$
is a lower semi-continuous proper convex function.
\item There exists a convex compact set $\Omega$, such that $F(t,x)\subseteq\Omega$
for all $(t,x)\in D$.%
\begin{comment}
The function $w(\cdot,\cdot)$ is either non-positive or non-negative
for all $(t,x)\in D$.
\end{comment}
\end{itemize}
\begin{comment}
%
\begin{rem}
\label{rem:sign w}The case with $w(\cdot,\cdot)\le0$ can be reduced
to the case with $w(t,x)\ge0$. Indeed, let $\nu(t,x)=\mu(1-t,x)$,
then 
\begin{align*}
\partial_{t}\nu(t,x)+\nabla\cdot\nu(t,x)(-v(1-t,x)) & =-w(1-t,x)\nu(t,x)
\end{align*}
Thus by introducing 
\[
\tilde{F}(t,x)=\left\{ (\alpha,\beta):\begin{bmatrix}-\alpha\\
-\beta
\end{bmatrix}\in F(1-t,x)\right\} 
\]
the problem can be reformulated as Problem \ref{pb:P} with non-negative
$w(\cdot,\cdot)$. Thus from now on, we assume that $w$ is non-negative. 
\end{rem}
\end{comment}

{\bfseries{}%
\begin{comment}
\textbf{Assumption H1}: The control input $u(t,x)$ takes values in
Let $v(t,x)=f(t,x,u(t,x))$ and $w(t,x)=g(t,x,u(t,x))$. Assume $v\in L^{1}([0,1];W^{1,\infty}(X;\mathbb{R}^{d}))$
and $w(t,\cdot)$ is bounded and locally Lipschitz for $t\in[0,1]$.

\textbf{Assumption H2}: the set $F(t,x)$ and the function $u\mapsto L(t,x,u)$
are both convex for each $(t,x)\in D$, and that there exists a compact
set $K$, such that $F(t,x)\subseteq K$ for every $(t,x)\in D$.
$g$ is bounded, $|g|\le C_{g}$.

Note that $M_{t}\le e^{C_{g}t}M_{0}$, where $M_{t}=\int_{X}{\rm d}\mu(t,x)$. 
\end{comment}
}

\begin{rem}
\label{rem:Fisher}Problem \ref{pb:P} extends the unbalanced optimal
transport distance proposed in \cite{chizatInterpolatingDistanceOptimal2018},
by taking $L(t,x,v,w)=\frac{1}{2}|v|^{2}+\frac{1}{2}\delta^{2}w^{2}$
and $F(t,x)=\mathbb{R}^{n+1}$.\footnote{Although $F(t,x)$ is non-compact, we can choose bounded regions which include the supports of $\mu_0, \mu_1$, and the evolution of the dynamics.}
In this case, the optimal cost is
denoted as $WF_{\delta}(\mu_{0},\mu_{1})^{2}$. The formulation (\ref{sys:mu_v_w})
is much general, which handles the dynamical constraints (\ref{eq:sys-cnst})
with more general cost (\ref{eq:cost-Lag}). In particular, system
(\ref{sys:fg}) with optimization cost (\ref{eq:cost}) can be turned
into Problem \ref{pb:P} by defining $L$ and $F$ as 
\begin{equation}
L(t,x,v,w):=\inf_{u}\left\{ L(t,x,u)\left|\begin{array}{c}
v=f(t,x,u)\\
w=g(t,x,u)
\end{array}\right.\right\} \label{eq:L->L}
\end{equation}
(we have used the same letter $L$ for two different functions but
it does not cause any confusion) and
\begin{equation}
F(t,x)=\left\{ \begin{bmatrix}v\\
w
\end{bmatrix}:v=f(t,x,u),\;w=g(t,x,u),\;\exists u\right\} .\label{eq:fg->F(t,x)}
\end{equation}
\end{rem}
Given the boundary condition $\mu_{0},\mu_{1}$, it is an essential
question whether Problem \ref{pb:P} is feasible. We have the following
sufficient condition.
\begin{lem}
\label{lem:control}Problem \ref{pb:P} is feasible for a given pair
$(\mu_{0},\mu_{1})$ if for every $(x_{0},z_{0})\in{\rm supp}\mu_{0}\times\mathbb{R}$
and $(x_{1},z_{1})\in{\rm supp}\mu_{1}\times\mathbb{R}$, there exists
a feasible pair $(v(\cdot,\cdot),w(\cdot,\cdot))$ such that $(x_{0},z_{0})$
and $(x_{1},z_{1})$ are the initial and terminal states of the system
$\dot{x}=v(t,x)$, $\dot{z}=w(t,x)$ respectively. 
In particular, if the system 
\begin{equation}
\left\{ \begin{array}{l}
\dot{x}=f(t,x,u)\\
\dot{z}=g(t,x,u)
\end{array}\right.\label{sys: x-m}
\end{equation}
is controllable, then by assuming that $F(t,x)$ is defined by (\ref{eq:fg->F(t,x)})
and $L$ as (\ref{eq:L->L}), the problem is feasible. 
\end{lem}
\section{Relaxation}

Problem \ref{pb:P} is an optimal control problem in measure space.
In classical optimal control theory in finite dimensional space, convex
duality theories have been systematically developed to turn the nonlinear
and non-convex problem into infinite dimensional linear programs (hence
convex), see for example, L. C. Young \cite{young1969lectures}, W.
Fleming and D. Vermes \cite{flemingConvexDualityApproach1989} and
R. Vinter \cite{vinterConvexDualityNonlinear1993}, see also A. Rantzer and P. A. Parrilo
\cite{rantzer2000convexity}. For balanced dynamical optimal transport,
some duality theories have also been established in the literature,
see for example \cite{bernardOptimalMassTransportation2007,ghoussoubOptimalTransportControlled2018}.
Duality results for unbalanced dynamical optimal transport
(UDOT) have been established by L. Chizat et. al. \cite{chizatInterpolatingDistanceOptimal2018,chizatUnbalancedOptimalTransport2017,chizatUnbalancedOptimalTransport2018} and M. Liero et. al. \cite{liero2018optimal}.
It was shown that the interpolating Wasserstein-Fisher-Rao
metric (see Remark \ref{rem:Fisher}) $WF_{\delta}$ can be formulated
as
\[
WF_{\delta}(\mu_{0},\mu_{1})^{2}=\sup_{\varphi\in C^{1}([0,1]\times X)}\int\varphi(1,\cdot){\rm d}\mu_{1}-\varphi(0,\cdot){\rm d}\mu_{0}
\]
where $\varphi$ is such that
\begin{equation}
\partial_{t}\varphi+\frac{1}{2}\left(|\partial_{x}\varphi|^{2}+\frac{\varphi^{2}}{\delta^{2}}\right)\le0.\label{eq:WF-HJB}
\end{equation}

This duality was obtained based on the fact that the control
$v$ and dissipation $w$ can be chosen freely. It is an open question
whether similar results hold when the control and dissipation are
constrained as in Problem \ref{pb:P}. The goal of this paper is to
derive the dual of Problem \ref{pb:P}. We demonstrate that, by leveraging
techniques from classical dual optimal control theory, a similar duality
result does hold. 

The key procedure in the construction of classical dual optimal control is to
lift the problem into a space composed by the so-called Young measures
\cite{young1969lectures}, or ``relaxed arcs'' \cite{vinterConvexDualityNonlinear1993}.
For example, in \cite{vinterConvexDualityNonlinear1993}, a relaxed
arc (when the terminal time is fixed) for the system $\dot{x}\in\Omega$
on $[0,1]$ with $x(0)=x_{0}$ is a tuple $(\gamma,m)$ where $\gamma:[0,1]\to X$
is Lipschitz on $[0,1]$ with $\gamma(0)=x_{0}$ and $m_{t}\in\mathscr{P}(X)$
is a one-parameter probability measure such that
\[
\dot{\gamma}(t)=\int_{\Omega}v{\rm d}m_{t}(v).
\]
To mimic the notion of relaxed arcs, we propose to define a notion
called {\it relaxed measures}. First, we consider a relaxation of Problem
\pageref{pb:P}:
\begin{problem}
\label{pb:relax-Omega}In Problem \ref{pb:P}, replace the constraint
(\ref{eq:sys-cnst}) by
\begin{equation}
\begin{bmatrix}v(t,x)\\
w(t,x)
\end{bmatrix}\in\Omega\label{eq:Omega-cnst}
\end{equation}
where $\Omega\subseteq\mathbb{R}^{n+1}$ is defined in Assumption H1.

We then define relaxation measures based on the formulation of Problem
\ref{pb:relax-Omega}.
\end{problem}
\begin{defn}[Relaxed measure]
\label{def:relax-meas}A pair $(\mu,m)$ is called a relaxed measure
if $\mu$ solves (\ref{sys:mu_v_w}) with $(v,w)$ being 
\[
\begin{bmatrix}v(t,x)\\
w(t,x)
\end{bmatrix}=\int_{\Omega}\begin{bmatrix}\alpha\\
\beta
\end{bmatrix}{\rm d}m_{t,x}(\alpha,\beta)
\]
where $m_{t,x}$ is a Young measure.%
\begin{comment}
supported on $F(t,x)$ for every $(t,x)\in[0,1]\times X$. 
\end{comment}
\end{defn}

A relaxed measure defines a bounded linear functional $\Phi$ on $C(D\times\Omega)$, for any $\xi \in C(D\times \Omega)$, define
\begin{equation}
\left\langle \Phi,\xi\right\rangle :=\int_{0}^{1}\int_{X}\int_{\Omega}\xi(t,x,\alpha,\beta){\rm d}m_{t,x}(\alpha,\beta){\rm d}\mu(t,x){\rm d}t.\label{eq:lin_func_Phi}
\end{equation}
\begin{comment}
Consider the system (\ref{sys:mu_v_w}), assume additionally that
the system (\ref{sys:mu_v_w}) satisfies
\[
|\mu|_{\mathcal{M}_{+}(D)}\le1,
\]
for all $(v,w)\in K$ and $\mu(0,\cdot)\in\Sigma_{0}\subseteq\mathcal{M}_{+}(X)$.
In this case, we say that $(\mu,m)$ is a relaxed arc for system (\ref{sys:mu_v_w}).
\end{comment}
\begin{comment}
\[
f(x,U)=\mathcal{F}(x)\subseteq\mathcal{F},\;g(x,U)=\mathcal{G}(x)\subseteq\mathcal{G},
\]
both $\mathcal{F}(x)$ and $\mathcal{G}(x)$ are bounded convex sets,
so are $\mathcal{F}$ and $\mathcal{G}$.
\end{comment}
The linear functional $\Phi$ has the following property: for every $\phi\in C^{1}([0,1]\times X)$, there holds
\begin{align}
\left\langle \Phi,\phi_{t}+\phi_{x}\cdot\alpha+\phi\beta\right\rangle  & =\int_{X}\phi(1,\cdot){\rm d}\mu_{1}-\phi(0,\cdot){\rm d}\mu_{0}.\label{eq:Phi_phi}
\end{align}
Indeed, 
\begin{align*}
 & \left\langle \Phi,\phi_{t}+\phi_{x}\cdot\alpha+\phi\beta\right\rangle \\
= & \int_{0}^{1}\int_{X}\int_{\Omega}\phi_{t}+\phi_{x}\cdot\alpha+\phi\beta{\rm d}m_{t,x}(\alpha,\beta){\rm d}\mu(t,x){\rm d}t\\
= & \int_{0}^{1}\int_{X}\left[\phi_{t}+\int_{\Omega}[\partial_{x}\phi,\;\phi]\cdot\begin{bmatrix}\alpha\\
\beta
\end{bmatrix}{\rm d}m_{t,x}(\alpha,\beta)\right]{\rm d}\mu(t,x){\rm d}t\\
= & \int_{0}^{1}\int_{X}\partial_{t}\phi+\partial_{x}\phi\cdot v+\phi w{\rm d}\mu(t,x){\rm d}t\\
= & \int_{X}\phi(1,\cdot){\rm d}\mu_{1}-\phi(0,\cdot){\rm d}\mu_{0}
\end{align*}

Conversely, we can also use \eqref{eq:Phi_phi} as a defining property of a set of linear functionals on $C(D\times\Omega)$, i.e., we set
\begin{equation}
W :=\{\Phi\in C^{*}(D\times\Omega)\;|\;\Phi\text{ satisfies (\ref{eq:Phi_phi})}\}\label{eq:set W}
\end{equation}
It is easy to observe that functionals in $W$ are non-negative in
the sense that $\left\langle \Phi,l\right\rangle \ge0$ when $0\le l\in C(D\times\Omega)$.
Now let 
\begin{equation}
S=\{\Phi\in C^{*}(D\times\Omega)\,|\,\exists\text{a relaxed measure s.t. (\ref{eq:lin_func_Phi}) holds}\}\label{eq:set S}
\end{equation}
then $S$ is a subset of $W$. The question is how big is the difference
between $S$ and $W$. The following lemma says that $W=\overline{{\rm co}}S$.
The proof is given in Appendix. 
\begin{lem}
\label{lem:W=00003Dco(S)}Let $W$ and $S$ be defined as (\ref{eq:set W})
and (\ref{eq:set S}) respectively. Then $W=\overline{{\rm co}}S$. 
\end{lem}
\begin{comment}

\subsection{Dirac mass}

Consider the system
\[
\partial_{t}\mu+\nabla\cdot\mu v=\mu w
\]
for $\mu_{0}=c_{0}\delta_{x_{0}}$. Consider the system
\[
\dot{x}=v(t,x)
\]
with $x(0)=x$. Let $x\mapsto X_{t}(x)$ be the system flow. Let 
\[
G_{t}(x)=\exp\left(\int_{0}^{t}w(s,X_{s}(x))ds\right)
\]
and 
\[
\rho_{t}(x)=c_{0}G_{t}(X_{t}(x))\delta_{X_{t}(x)}.
\]
Let's verify that $\rho_{t}(x)$ is the solution.

First, $\rho_{0}(x)=c_{0}\delta_{0}=\mu_{0}$. Next, 
\begin{align*}
\int_{0}^{t}\int_{X}\partial_{t}\varphi+\nabla\varphi\cdot v+\varphi wd\rho_{s}(x) & =\int_{0}^{t}c_{0}G_{s}(X_{s}(x))[\partial_{t}\varphi+\nabla\varphi\cdot v+\varphi w]ds\\
 & =\int_{0}^{t}c_{0}\frac{d}{ds}G_{s}(X_{s}(x))\varphi(s,X_{s}(x))ds\\
 & =c_{0}(G_{t}(X_{t})\varphi(t,X_{t}(x))-\varphi(0,x))\\
 & =\int\varphi(t,x)d\rho_{t}(x)-\int\varphi(0,x)d\rho_{0}(x)
\end{align*}

Let $m$ be a Young measure. 
\end{comment}

With Lemma \ref{lem:W=00003Dco(S)} at hand, we propose to study:
\begin{problem}[Relaxed primal problem]
\label{pb:relax}
\[
\min_{\Phi\in W}\left\langle \Phi,L\right\rangle 
\]
with constraint\footnote{Remember that $d_{F(t,x)}(\alpha, \beta)$ is the distance from $(\alpha, \beta)$ to the set $F(t,x)$, see the {\it Notation} paragraph of Section \ref{sec:intro}.}
\begin{align}
\left\langle \Phi,d_{F(t,x)}(\alpha,\beta)\right\rangle  & =0.\label{eq:Phi, d_F}
\end{align}
\end{problem}
The following lemma justifies that Problem \ref{pb:relax} is a tight
relaxation of Problem \ref{pb:P}. 
\begin{lem}
\label{lem:P=00003DW}The optimal values of Problem \ref{pb:P} and
Problem \ref{pb:relax} are equal.
\end{lem}

The complete proof of the lemma is tedious; a proof sketch
is provided in the Appendix, which is based on Lemma \ref{lem:W=00003Dco(S)}

\section{Duality}

In this section, we work toward a convex duality of Problem \ref{pb:relax}. We shall prove that
Problem \ref{pb:relax} admits the following dual formulation.
\begin{problem}[Dual problem]
\label{pb:Dual}
\[
\sup_{\phi}\int_{X}\phi(1,\cdot){\rm d}\mu_{1}-\phi(0,\cdot){\rm d}\mu_{0}
\]
where $\phi$ is taken over the $C^{1}$ solutions to the Hamilton-Jacobi-Bellman
inequality 
\begin{equation}
\partial_{t}\phi+H(t,x,\partial_{x}\phi,\phi)\le0\label{eq:sub HJB}
\end{equation}
with $H$ being defined as
\[
H(t,x,p,p_{0}):=\sup_{(v,w)\in F(t,x)}\{p^{\top}v+p_{0}w-L(t,x,v,w)\}.
\]
\end{problem}
Below the main result of this paper, whose proof is deferred to the Appendix:
\begin{thm}
\label{thm:main}The optimal values of Problems \ref{pb:P}, \ref{pb:relax}
and \ref{pb:Dual} are equal.
\end{thm}

Now the dual formulation of $WF_{\delta}$ (see \eqref{eq:WF-HJB}) is now justified by Theorem \ref{thm:main}. The Hamiltonian reads
\begin{align*}
H(t,x,p,p_{0}) & =\sup_{(v,w)}\left\{ p^{\top}v+p_{0}w-\left(\frac{1}{2}|v|^{2}+\frac{1}{2}\delta^{2}w^{2}\right)\right\} \\
 & =\frac{1}{2}\left(|p|^{2}+\frac{p_{0}^{2}}{\delta^{2}}\right)
\end{align*}
based on which we immediately recover (\ref{eq:WF-HJB}) from (\ref{eq:sub HJB}).

Theorem \ref{thm:main} also recovers the well-known duality result
on balanced optimal transport, i.e., when $w\equiv0$, first proved by Bernard and 
Buffoni \cite{bernardOptimalMassTransportation2007}. We have a more generalized version:
\begin{cor}
Consider the differential inclusion
\[
\dot{x}\in G(t,x)
\]
on $X$, with $G(t,x)\subseteq\Omega$ being convex and $\Omega$
some convex compact set. Let 
\[
c(x,y)=\inf\int_{0}^{1}L(t,x,\dot{x}){\rm d}t
\]
for some Lagrangian $L$ such that $v\mapsto L(t,x,v)$ is convex.
Let $\mu_{0},\mu_{1}$ be two probability measures, then 
\[
\min_{\pi\in\Gamma(\mu_{0},\mu_{1})}\int c(x,y){\rm d}\pi(x,y)=\sup_{\phi}\int\phi_{1}(1\cdot){\rm d}\mu_{1}-\phi_{0}(0,\cdot){\rm d}\mu_{0}
\]
where the supremum is taken over the set of $C^{1}$ solutions to
the HJB inequality 
\[
\partial_{t}\phi+H(t,x,\partial_{x}\phi)\le0.
\]
\end{cor}

\section{A Numerical Algorithm}

We propose a numerical algorithm similar to the one proposed in \cite{benamouComputationalFluidMechanics2000}
based on the dual formulation Problem \ref{pb:Dual}. To streamline
the exposition, we introduce some notations below
\begin{align*}
G(\phi) & =\int_{X}\phi(0,\cdot){\rm d}\mu_{0}-\phi(1,\cdot){\rm d}\mu_{1}\\
K_{H} & \;=\{(a,b,c)\in C(D;\mathbb{R}^{n+2}):a+H(t,x,b,c)\le0\}
\end{align*}
Then, Problem \ref{pb:Dual} can be written as
\[
\sup_{q,\phi}-G(\phi)-I_{K_{H}}(q)
\]
subject to the constraint $q=(\partial_{t}\phi,\partial_{x}\phi,\phi)$.
We introduce a Lagrangian multiplier $\lambda:D\to\mathbb{R}^{n+2}$
to accommodate this equality constraint:
\[
\sup_{q,\phi}\inf_{\lambda}-G(\phi)-I_{K_{H}}(q)-\left\langle \lambda,q-\delta(\phi)\right\rangle 
\]
where $\delta$ is the operator 
\[
\delta(\phi)=\begin{bmatrix}\partial_{t}\phi\\
\partial_{x}\phi\\
\phi
\end{bmatrix}.
\]
Let 
\[
\mathcal{L}(\phi,q,\lambda):=G(\phi)+I_{K_{H}}(q)+\left\langle \lambda,\delta(\phi)-q\right\rangle 
\]
then we shall study $\inf_{q,\phi}\sup_{\lambda}\mathcal{L}(\phi,q,\lambda)$
or its augmented Lagrangian form
\[
\inf_{q,\phi}\sup_{\lambda}\mathcal{L}_{r}(\phi,q,\lambda)
\]
where $\mathcal{L}_{r}(\phi,q,\lambda)=\mathcal{L}(\phi,q,\lambda)+\frac{r}{2}\left\Vert \delta(\phi)-q\right\Vert ^{2}$.
The norm $\|\cdot\|$ is defined as $\|v\|^{2}=\int_{D}|v(t,x)|^{2}{\rm d}t{\rm d}x$. 

\begin{algorithm}[h]
Given $(\phi^{0},q^{0},\lambda^{1})$, iterate for $k\ge1$
\begin{enumerate}
\item Compute $\phi^{k}\in\arg\min_{\phi}\mathcal{L}_{r}(\phi,q^{k-1},\lambda^{k})$;
\item Compute $q^{k}\in\arg\min_{q}\mathcal{L}_{r}(\phi^{k},q,\lambda^{k})$;
\item $\lambda^{k+1}=\lambda^{k}+r(\delta(\phi^{k})-q^{k})$.
\end{enumerate}
\caption{Augmented Lagrangian Method\label{alg:ALM}}
\end{algorithm}

In Step 1), $\phi^{k}$ can be obtained by solving the second order
PDE
\[
-r(\Delta_{t,x}\phi-\phi)=\xi^{k}
\]
where
\begin{align*}
\xi^{k} & =-r(\partial_{t}a^{k-1}+\nabla_{x}\cdot b^{k-1}-c^{k-1})\\
 & \quad\quad+\nabla_{x}\cdot\lambda_{b}^{k}+\partial_{t}\lambda_{a}^{k}-\partial_{t}\lambda_{c}^{k}
\end{align*}
with pure Neumann boundary conditions
\begin{align*}
r\partial_{t}\phi(1,\cdot) & =ra^{k-1}(1,\cdot)+\mu_{1}-\lambda_{a}^{k}(1,\cdot)\\
r\partial_{t}\phi(0,\cdot) & =ra^{k-1}(0,\cdot)+\mu_{0}-\lambda_{a}^{k}(0,\cdot)
\end{align*}
and $\partial_{x}\phi=0$ on $\partial X$. Note we have denoted $\lambda^{k}=(\lambda_{a}^{k},\lambda_{b}^{k},\lambda_{c}^{k})$
and $q^{k}=(a^{k},b^{k},c^{k})$.
To obtain this, differentiating $\mathcal{L}_{r}$ w.r.t. $\phi$:
\begin{equation}
G(\varphi)+\left\langle \lambda^{k},\delta(\varphi)\right\rangle +r\left\langle \delta(\phi)-q^{k-1},\delta(\varphi)\right\rangle =0\label{eq:G(varphi)}
\end{equation}
which should hold for all $\varphi$.%
\begin{comment}
\begin{align*}
 & \left\langle \delta(\phi)-q^{n-1},\delta(\varphi)\right\rangle \\
= & \int(\phi_{t}-a^{n-1})\varphi_{t}+(\phi_{x}-b^{n-1})\cdot\varphi_{x}+(\phi-c^{n-1})\varphi{\rm d}t{\rm d}x\\
=- & \int(\phi_{tt}-a_{t}^{n-1})\varphi{\rm d}t{\rm d}x+\int\varphi(1,\cdot)(\phi_{t}(1,\cdot)-a^{n-1}(1,\cdot)){\rm d}x\\
 & -\int\varphi(0,\cdot)(\phi_{t}(0,\cdot)-a^{n-1}(0,\cdot)){\rm d}x\\
 & -\int(\phi_{xx}-\nabla_{x}\cdot b^{n-1})\varphi{\rm d}t{\rm d}x+\int_{\partial X}\varphi(\phi_{x}-b^{n-1})\cdot dS\\
 & +\int(\phi-c^{n-1})\varphi{\rm d}t{\rm d}x.
\end{align*}
\end{comment}
{} Easy calculations show that
\begin{align*}
\left\langle \lambda^{k},\delta(\varphi)\right\rangle = & \int\lambda_{a}^{k}\varphi_{t}+\lambda_{b}^{k}\cdot\varphi_{x}+\lambda_{c}^{k}\varphi{\rm d}x{\rm d}t\\
= & -\int(\lambda_{a}^{k})_{t}\varphi{\rm d}x{\rm d}t\\
 & +\int\lambda_{a}^{k}(1,\cdot)\varphi(1,\cdot)-\lambda_{a}^{k}(0,\cdot)\varphi(0,\cdot){\rm d}x\\
 & -\int(\nabla_{x}\lambda_{b}^{k})\varphi{\rm d}x{\rm d}t+\int_{\partial X}\varphi\lambda_{b}^{k}\cdot dS\\
 & +\int\lambda_{c}^{k}\varphi{\rm d}x{\rm d}t
\end{align*}
and
\begin{align*}
 & \left\langle \delta(\phi)-q^{k-1},\delta(\varphi)\right\rangle \\
= & -\int r(\phi_{tt}+\phi_{xx}-\phi-a_{t}^{k-1}-\nabla_{x}\cdot b^{k-1}+c^{k-1})\varphi{\rm d}t{\rm d}x\\
 & +\int\varphi(1,\cdot)(\phi_{t}(1,\cdot)-a^{k-1}(1,\cdot)){\rm d}x\\
 & -\int\varphi(0,\cdot)(\phi_{t}(0,\cdot)-a^{k-1}(0,\cdot)){\rm d}x\\
 & -\int_{\partial X}(\phi_{x}-b^{k-1})\varphi\cdot dS
\end{align*}
Substituting these into (\ref{eq:G(varphi)}) gives the desired result.
\begin{comment}
\[
-r(\Delta_{t,x}\phi-\phi)+r(\partial_{t}a^{k-1}+\nabla_{x}\cdot b^{k-1}-c^{k-1})-\nabla_{x}\cdot\lambda_{b}^{k}-\partial_{t}\lambda_{a}^{k}+\partial_{t}\lambda_{c}^{k}=0
\]
\end{comment}
\begin{comment}
\[
-r(\Delta_{t,x}\phi-\phi)=\xi_{k}
\]
where
\begin{align*}
\xi_{k} & =-r(\partial_{t}a^{k-1}+\nabla_{x}\cdot b^{k-1}-c^{k-1})\\
 & \quad\quad+\nabla_{x}\cdot\lambda_{b}^{k}+\partial_{t}\lambda_{a}^{k}-\partial_{t}\lambda_{c}^{k}
\end{align*}
and 
\begin{align*}
r\partial_{t}\phi(1,\cdot) & =ra^{k-1}(1,\cdot)+\rho_{1}-\lambda_{a}^{k}(1,\cdot)\\
r\partial_{t}\phi(0,\cdot) & =ra^{k-1}(0,\cdot)+\rho_{0}-\lambda_{a}^{k}(0,\cdot)
\end{align*}
and on $\partial X$:
\[
\partial_{x}\phi=0
\]
by assuming that $\mu$ is compactly supported in $X$. 
\end{comment}

Step 2) amounts to the pointwise minimization
\begin{align*}
q^{k} & =\arg\inf_{q\in K_{H}}\frac{r}{2}|q-\delta(\phi^{k})|^{2}+\left\langle \lambda^{k},\delta(\phi^{k})-q\right\rangle \\
 & =\arg\inf_{q\in K_{H}}\frac{r}{2}|q|^{2}-\left\langle q,\lambda^{k}+r\delta(\phi^{k})\right\rangle \\
 & =\arg\inf_{q\in K_{H}}\left|q-(\lambda^{k}/r+\delta(\phi^{k}))\right|^{2}.
\end{align*}

\section{Conclusion}

We have proposed a generalized dynamical unbalanced optimal transport (UOT)
framework with constrained control and mass changes. By deriving its
convex duality using classical dual optimal control techniques, a
first-order numerical algorithm has been proposed. Our approach enhances
UOT\textquoteright s applicability to control problems like energy
systems and biological modeling. Future work includes proposing more
efficient numerical algorithms to handle the dimensionality issue
associated with the solution of the HJB.

\section*{Appendix}
\begin{proof}[Proof of Lemma \ref{lem:W=00003Dco(S)}]
Without loss of generality, assume $\int{\rm d}\mu_{0}=1$. Suppose
that there exists $\Lambda_{0}\in\mathscr{W}\backslash\overline{{\rm co}}S$,
then there exists $l\in C(D\times\Omega)$ and $a\in\mathbb{R}$,
such that
\begin{equation}
\left\langle \Phi,l\right\rangle \ge a,\quad\forall\Phi\in S\label{eq:Phi-l a}
\end{equation}
and
\[
\left\langle \Lambda_{0},l\right\rangle <a.
\]
\begin{comment}
If we can show
\begin{align}
a & \le\int\phi_{0}{\rm d}\mu_{0}-\phi_{1}{\rm d}\mu_{1}\label{eq:a phi0 phi1}
\end{align}
for some $\phi\in C^{1}$ satisfying
\[
-\phi_{t}-\phi_{x}\cdot\alpha-\phi\beta-l(t,x,\alpha,\beta)\le0,\quad\forall\alpha,\beta
\]
then because
\[
a>\left\langle \Lambda_{0},l\right\rangle \ge-\left\langle \Lambda_{0},\phi_{t}+\phi_{x}\cdot\alpha+\phi\beta\right\rangle =\int\phi_{0}{\rm d}\mu_{0}-\phi_{1}{\rm d}\mu_{1}\ge a
\]
we would be able to derive a contradiction. 
\end{comment}

Given $x\in X$, define 
\[
\phi(t,x)=\inf_{\mu,m}\int_{t}^{1}\int_{X}\int_{\Omega}l(s,z,\alpha,\beta){\rm d}m_{s,z}(\alpha,\beta){\rm d}\mu(s,z){\rm d}s
\]
with $\mu(t,\cdot)=\delta_{x}$ and $(\mu,m)$ is a relaxed measure.
By definition, $\phi(0,x)\ge a$, and thus 
\begin{equation}
a\le\int\phi(0,\cdot){\rm d}\mu_{0}-\phi(1,\cdot){\rm d}\mu_{1}\label{eq:a =00005Cleq}
\end{equation}
since $\mu_{0}$ is a probability measure and $\phi(1,\cdot)=0$.

For fixed $(\alpha,\beta)\in\Omega$, let $m_{s,x}=\delta_{(\alpha,\beta)}$,
$\forall s\in[t,1]$. Then $x(s)=x+(s-t)\alpha$, and $\mu(s,\cdot)=p(s)\delta_{x(s)}$
with $p(s)=e^{\beta(s-t)}$. By Bellman's principle of optimality,
{\small
\[
\phi(t,x)\le p(t+h)\phi(t+h,x(t+h))+\int_{t}^{t+h}l(s,x(s),\alpha,\beta){\rm d}s
\]
}or{\small
\[
\frac{\phi(t,x)-p(t+h)\phi(t+h,x(t+h))}{h}\le\frac{1}{h}\int_{t}^{t+h}l(s,x(s),\alpha,\beta){\rm d}s
\]
}If $\phi$ is $C^{1}$, we deduce 
\[
-\partial_{t}\phi(t,x)-\partial_{x}\phi(t,x)\cdot\alpha-\phi(t,x)\beta\le l(t,x,\alpha,\beta)
\]
by letting $h\to0+$. Since $t,x,\alpha,\beta$ can be chosen freely,
this inequality holds for all $(t,x,\alpha,\beta)$. Now 
\begin{align*}
a>\left\langle \Lambda_{0},l\right\rangle  & \ge-\left\langle \Lambda_{0},\partial_{t}\phi+\partial_{x}\phi\cdot\alpha-\phi\beta\right\rangle \\
 & =\int\phi(0,\cdot){\rm d}\mu_{0}-\phi(1,\cdot){\rm d}\mu_{1}
\end{align*}
which yields a contradiction by (\ref{eq:a =00005Cleq}). To complete
the proof, we approximate $\phi$ by $C^{1}$ functions in a similar
fashion as in \cite{vinterConvexDualityNonlinear1993}.
\end{proof}
\begin{proof}[Proof sketch of Lemma \ref{lem:P=00003DW}]
Given $\Phi\in W$, express $\Phi$ as 
\[
\Phi=\int_{S}\Lambda{\rm d}P(\Lambda)
\]
for some Borel probability measure $P$ on $S$. For any $l$, $\left\langle \Phi,l\right\rangle =\int_{S}\left\langle \Lambda,l\right\rangle {\rm d}P(\Lambda)$
implies the existence of some $\Lambda_{0}\in S$ such that $\left\langle \Phi,l\right\rangle =\left\langle \Lambda_{0},l\right\rangle =\int l{\rm d}m_{t,x}{\rm d}\mu{\rm d}t$.
By (\ref{eq:Phi, d_F}), we have
\[
\int_{S}\left\langle \Lambda,d_{F(t,x)}(\alpha,\beta)\right\rangle {\rm d}P(\Lambda)=\left\langle \Phi,d_{F(t,x)}(\alpha,\beta)\right\rangle =0
\]
from which we deduce $\left\langle \Lambda,d_{F(t,x)}(\alpha,\beta)\right\rangle =0$
a.e. $\Lambda$. Therefore, we may assume $\left\langle \Lambda_{0},d_{F(t,x)}(\alpha,\beta)\right\rangle =0$,
implying that $m_{t,x}$ is supported in $F(t,x)$ for ${\rm d}\mu{\rm d}t$
a.e. $(t,x)$. Then 
\[
\begin{bmatrix}v(t,x)\\
w(t,x)
\end{bmatrix}:=\int_{\Omega}\begin{bmatrix}\alpha\\
\beta
\end{bmatrix}{\rm d}m_{t,x}(\alpha,\beta)\in F(t,x),\;\text{a.e. }(t,x)
\]
and 
\[
\int_{\Omega}L(t,x,\alpha,\beta){\rm d}m_{t,x}(\alpha,\beta)\ge L(t,x,v(t,x),w(t,x))
\]
by Jensen's inequality. To summarize, we have constructed from an
arbitrary $\Phi\in W$ under constraint (\ref{eq:Phi, d_F}) a tuple
$(\mu,v,w)$ satisfying the equation (\ref{sys:mu_v_w}) in a way
that 
\[
\left\langle \Phi,L\right\rangle \ge\int_{0}^{1}\int_{X}L(t,x,v,w){\rm d}\mu(t,x){\rm d}t,
\]
implying that the optimal value of Problem \ref{pb:P} is upper bounded
by that of Problem \ref{pb:relax}. On the other hand, the reverse
inequality holds in an obvious way. This gives the desired result.

\end{proof}
\begin{proof}[Proof of Theorem \ref{thm:main}]
First, assume that $w(t,x)\ge0$ for all $(v(t,x),w(t,x))\in F(t,x)$.
The assumption will be removed later. Consider the space $\mathcal{X}=C(D\times\Omega)\times C(X)\times C(X)$
and its topological dual $\mathcal{X}^{*}=C^{*}(D\times\Omega)\times C^{*}(X)\times C^{*}(X)$.
We construct two convex functionals on $\mathcal{X}^{*}$ as follows.
First, define two subsets in $\mathcal{X}^{*}$:
\begin{align*}
P & =\left\{ (m,n_{0},n_{1})\in\mathcal{X}^{*}\left|\begin{array}{c}
n_{0}=\mu_{0},\;n_{1}=\mu_{1},\;m\ge0\\
\left\langle m,d_{F(t,x)}(\cdot,\cdot)\right\rangle =0
\end{array}\right.\right\} \\
Q & =\left\{ (m,n_{0},n_{1})\in\mathcal{X}^{*}\left|\begin{array}{c}
\text{(\ref{eq:m_n1_n0}) holds for all}\\
\phi\in C^{1},\,n_{0},\,n_{1}\in C^{*}(X)
\end{array}\right.\right\} 
\end{align*}
\begin{equation}
\left\langle m,\partial_{t}\phi+\phi_{x}\cdot v+\phi w\right\rangle =\left\langle n_{1},\phi_{1}\right\rangle -\left\langle n_{0},\phi_{0}\right\rangle \label{eq:m_n1_n0}
\end{equation}
Then define two functionals $p^{*}$, $q^{*}$ (the asterisk superscripts
is to remind that they are defined on $\mathcal{X}^{*}$) as
\begin{align*}
p^{*}(m,n_{0},n_{1}) & =\left\langle m,L\right\rangle +I_{P}(m,n_{0},n_{1}),\\
q^{*}(m,n_{0},n_{1}) & =I_{Q}(m,n_{0},n_{1})
\end{align*}
Note that $q^{*}(m,n_{0},n_{1})=q^{*}(-m,-n_{0},-n_{1})$, it is then
easy to see that Problem \ref{pb:relax} is equivalent to 
\[
\min_{(m,n_{0},n_{1})\in\mathcal{X}^{*}}p^{*}(m,n_{0},n_{1})+q^{*}(-m,-n_{0},-n_{1})
\]
We shall use the Fenchel-Rockafellar theorem (Theorem \ref{thm:Fen-Rock}
in Appendix) to derive the dual of the above problem. First, we compute
$p$ and $q$ such that $p^{*}$ and $q^{*}$ are their convex conjugates. 

The computation of $p$ is rather straightforward, 
\begin{align*}
 & p(\xi,\eta_{0},\eta_{1})\\
= & \sup_{(m,n_{0},n_{1})}\left\langle m,\xi\right\rangle +\left\langle n_{0},\eta_{0}\right\rangle +\left\langle n_{1},\eta_{1}\right\rangle -p^{*}(m,n_{0},n_{1})\\
= & \sup_{m\ge0,\,\left\langle m,d_{F(t,x)}\right\rangle =0}\left\langle m,\xi-L\right\rangle +\int\eta_{0}{\rm d}\mu_{0}+n_{1}{\rm d}\mu_{1}\\
= & \max_{(t,x)\in D,\;(v,w)\in F(t,x)}(\xi-L)^{+}+\int\eta_{0}{\rm d}\mu_{0}+n_{1}{\rm d}\mu_{1}
\end{align*}
where $(\xi-L)^{+}:=\max\{0,\xi-L\}$. To compute $q$, consider the
set 
\[
\mathcal{Q}=\left\{ (\xi,\eta_{0},\eta_{1})\in\mathcal{X}\left|\begin{array}{c}
\exists\phi\in C^{1},\text{ s.t. }\eta_{0}=\phi_{0},\eta_{1}=-\phi_{1}\\
\xi=\phi_{t}+\phi_{x}\cdot\alpha+\phi\beta,\;\forall(t,x,\alpha,\beta)
\end{array}\right.\right\} 
\]
We assert that 
\[
q(\xi,\eta_{0},\eta_{1})=I_{\bar{\mathcal{Q}}}
\]
with $\bar{\mathcal{Q}}$ being the closure of $Q$. In fact, by definition
\[
q(\xi,\eta_{0},\eta_{1})=\sup_{(m,n_{0},n_{1})\in Q}\left\langle m,\xi\right\rangle +\left\langle \eta_{0},n_{0}\right\rangle +\left\langle \eta_{1},n_{1}\right\rangle 
\]
if $(\xi,\eta_{0},\eta_{1})\in\mathcal{Q}$, we have 
\begin{align*}
\left\langle m,\xi\right\rangle  & =\left\langle m,\phi_{t}+\phi_{x}\cdot\alpha+\phi\beta\right\rangle \\
 & =\left\langle n_{1},\phi_{1}\right\rangle -\left\langle n_{0},\phi_{0}\right\rangle \text{ (by (\ref{eq:m_n1_n0}))}\\
 & =-\left\langle n_{1},\eta_{1}\right\rangle -\left\langle n_{0},\eta_{0}\right\rangle \text{ (by definition of \ensuremath{\mathcal{Q}})}
\end{align*}
from which it follows that $q(\xi,\eta_{0},\eta_{1})=0$. If $(\xi,\eta_{0},\eta_{1})\notin\bar{\mathcal{Q}}$,
then since $\mathcal{Q}$ is a closed linear space, there exists $(m',n_{0}',n_{1}')\in\mathcal{X}^{*}$
such that 
\begin{equation}
\left\langle \xi,m'\right\rangle +\left\langle \eta_{0},n_{0}'\right\rangle +\left\langle \eta_{1},n_{1}'\right\rangle >0\label{eq:sep}
\end{equation}
and
\[
\left\langle \partial_{t}\phi+\phi_{x}\cdot v+\phi w,m'\right\rangle +\left\langle \phi_{0},n_{0}'\right\rangle +\left\langle -\phi_{1},n_{1}'\right\rangle =0
\]
for all $\phi\in C^{1}$. But this implies $(m',n_{0}',n_{1}')\in Q$.
By (\ref{eq:sep}), we arrive at the conclusion that $q(\xi,\eta_{0},\eta_{1})=+\infty$
since $\alpha(m',n_{0}',n_{1}')\in Q$ for any $\alpha\in\mathbb{R}$. 

Now invoking Theorem \ref{thm:Fen-Rock}, we have
\begin{align}
 & \min_{(m,n_{0},n_{1})\in\mathcal{X}^{*}}p^{*}(m,n_{0},n_{1})+q^{*}(-m,-n_{0},-n_{1})\nonumber \\
= & \sup_{(\xi,\eta_{0},\eta_{1})\in\mathcal{X}}-p(\xi,\eta_{0},\eta_{1})-q(\xi,\eta_{0},\eta_{1})\nonumber \\
= & \sup_{(\xi,\eta_{0},\eta_{1})\in\bar{\mathcal{Q}}}-p(\xi,\eta_{0},\eta_{1})\nonumber \\
= & \sup_{(\xi,\eta_{0},\eta_{1})\in\bar{\mathcal{Q}}}\left\{ -\max(\xi-L)^{+}-\int\eta_{0}{\rm d}\mu_{0}-n_{1}{\rm d}\mu_{1}\right\} \nonumber \\
= & \sup_{\phi\in C^{1}}-\max(\phi_{t}+\phi_{x}\cdot\alpha+\phi\beta-L)^{+}+\int\phi_{1}{\rm d}\mu_{1}-\phi_{0}{\rm d}\mu_{0}\label{eq:sup -max}
\end{align}
where the maximization is taken over the set $(t,x,\alpha,\beta)\in D\times\Omega$
such that $(\alpha,\beta)\in F(t,x)$.

For fixed $\phi\in C^{1}$, let $c=\max(\phi_{t}+\phi_{x}\cdot\alpha+\phi\beta-L)^{+}$
and define $\tilde{\phi}=\phi-ct$. Then
\begin{align*}
 & \tilde{\phi}_{t}+\tilde{\phi}_{x}\cdot\alpha+\tilde{\phi}\beta-L\\
= & (\phi_{t}+\phi_{x}\cdot\alpha+\phi\beta-L-c)-ct\beta\\
\le & 0
\end{align*}
since $\beta\ge0$ by assumption. Note that this is equivalent to
saying that 
\begin{equation}
\partial_{t}\tilde{\phi}+H(t,x,\partial_{x}\tilde{\phi},\tilde{\phi})\le0.\label{eq:tilde phi HJB}
\end{equation}
By assuming without loss of generality that $\mu_{1}(X)=1$, we get
\begin{align}
 & -\max(\tilde{\phi}_{t}+\tilde{\phi}_{x}\cdot\alpha+\tilde{\phi}\beta-L)^{+}+\int\tilde{\phi}_{1}{\rm d}\mu_{1}-\tilde{\phi}_{0}{\rm d}\mu_{0}\nonumber \\
= & \int\tilde{\phi}_{1}{\rm d}\mu_{1}-\tilde{\phi}_{0}{\rm d}\mu_{0}\label{eq:tilde cost}\\
= & \int(\phi_{1}-c){\rm d}\mu_{1}-\phi_{0}{\rm d}\mu_{0}\nonumber \\
= & -c+\int\phi_{1}{\rm d}\mu_{1}-\phi_{0}{\rm d}\mu_{0}\nonumber \\
= & -\max(\phi_{t}+\phi_{x}\cdot\alpha+\phi\beta-L)^{+}+\int\phi_{1}{\rm d}\mu_{1}-\phi_{0}{\rm d}\mu_{0}\nonumber 
\end{align}
Now, replacing $\phi$ in (\ref{eq:sup -max}) by $\tilde{\phi}$
and in view of (\ref{eq:tilde phi HJB}), (\ref{eq:tilde cost}),
the theorem is proven under the assumption that $w(\cdot,\cdot)\ge0$. 

To remove the positivity assumption on $w$, consider the scaling
$\nu(t,x)=e^{\lambda t}\mu(t,x)$. Then $\nu$ satisfies
\begin{align*}
\partial_{t}\nu+\nabla\cdot\nu v & =\tilde{w}\nu\\
\nu(0,\cdot) & =\mu_{0}\\
\nu(1,\cdot) & =e^{\lambda}\mu_{1}
\end{align*}
in which $\tilde{w}=w+\lambda$. The minimization of (\ref{eq:cost-Lag})
becomes 
\[
\int_{0}^{1}\int_{X}L(t,x,v,\tilde{w}-\lambda)e^{-\lambda t}{\rm d}\nu(t,x){\rm d}t.
\]
Choose $\lambda$ large enough such that $\tilde{w}\ge0$, which is
possible due to the compactness of $\Omega$. Applying the previous
result, we get a dual formulation of this problem:
\begin{equation}
\sup_{\phi}\int_{X}\phi(1,\cdot)e^{\lambda}{\rm d}\mu_{1}-\phi(0,\cdot){\rm d}\mu_{0}\label{eq:sup phi lambda}
\end{equation}
subject to
\begin{equation}
\partial_{t}\phi+\tilde{H}(t,x,\partial_{x}\phi,\phi)\le0\label{eq:HJB tilde}
\end{equation}
The Hamiltonian $\tilde{H}$ is computed as follows
\begin{align*}
 & \tilde{H}(t,x,p,p_{0})\\
= & \sup_{(v,\tilde{w})}p^{\top}v+p_{0}\tilde{w}-e^{-\lambda t}L(t,x,v,\tilde{w}-\lambda)\\
= & p_{0}\lambda+\sup_{(v,w)\in F(t,x)}p^{\top}v+p_{0}w-e^{-\lambda t}L(t,x,v,w)\\
= & p_{0}\lambda+e^{-\lambda t}\sup_{(v,w)\in F(t,x)}(e^{\lambda t}p)^{\top}v+(e^{\lambda t}p_{0})w-L\\
= & p_{0}\lambda+e^{-\lambda t}H(t,x,e^{\lambda t}p,e^{\lambda t}p_{0})
\end{align*}
where $H$ is the Hamiltonian defined as in Problem \ref{pb:Dual}.
Thus the HJB inequality can be transformed into 
\[
\partial_{t}\phi+\lambda\phi+e^{-\lambda t}H(t,x,e^{\lambda t}\partial_{x}\phi,e^{\lambda t}\phi)\le0
\]
or
\begin{equation}
\partial_{t}\varphi+H(t,x,\partial_{x}\varphi,\varphi)\le0\label{eq:HJB varphi}
\end{equation}
with $\varphi=e^{\lambda t}\phi$. Now, replacing $\phi$ in (\ref{eq:sup phi lambda})
by $e^{-\lambda t}\varphi$, (\ref{eq:sup phi lambda}) can be written
as
\[
\sup_{\varphi\in C^{1}}\int\varphi(1,\cdot){\rm d}\mu_{1}-\varphi(0,\cdot){\rm d}\mu_{0}
\]
in which $\varphi$ is subject to the HJB inequality (\ref{eq:HJB varphi}).
The proof is now complete.

\begin{comment}
\begin{align*}
 & \min_{(m,n_{0},n_{1})}p(m,n_{0},n_{1})+q(m,n_{0},n_{1})\\
= & \max_{(\xi,\eta_{0},\eta_{1})}-p^{*}(\xi,\eta_{0},\eta_{1})-q^{*}(-\xi,-\eta_{0},-\eta_{1})\\
= & \max_{\phi}\left\{ -\max(\xi-L)^{+}-\int\eta_{0}{\rm d}\mu_{0}-\int\eta_{1}{\rm d}\mu_{1}\right\} \\
= & \max_{\varphi}\left\{ -\max(\partial_{t}\varphi+\varphi_{x}\cdot\alpha+\varphi\beta-L)^{+}+\int\varphi_{0}{\rm d}\mu_{0}-\int\varphi_{1}{\rm d}\mu_{1}\right\} \\
= & \max_{\varphi}\int\varphi_{0}{\rm d}\mu_{0}-\int\varphi_{1}{\rm d}\mu_{1}
\end{align*}
\end{comment}
\end{proof}
\begin{thm}[Fenchel-Rockafellar]
\label{thm:Fen-Rock}Let $\mathcal{X}$ be a normed vector space,
$\mathcal{X}^{*}$ its topological dual. Let $p,q:\mathcal{X}\to\mathbb{R}\cup\{+\infty\}$
be two proper convex functionals. Denote $p^{*},q^{*}$ the convex
conjugates of $p,q$ respectively. If there exists a point $x_{0}\in\mathcal{X}$,
such that $p$ is continuous in at $x_{0}$ and $q(x_{0})<+\infty$,
then 
\[
\min_{x^{*}\in\mathcal{X}^{*}}\{p^{*}(x^{*})+q^{*}(-x^{*})\}=\sup_{x\in\mathcal{\mathcal{X}}}\{-p(x)-q(x)\}.
\]
\end{thm}

\section*{Acknowledgment}
The first author gratefully acknowledges the insightful discussions with Richard Vinter during a visit to Imperial College London, which significantly shaped this work.

\bibliographystyle{IEEEtran}
\bibliography{IEEEabrv,partialOT}

\end{document}